\documentclass[11pt]{article}

\newfont{\bb}{msbm10 at 11pt}
\def\r{\hbox{\bb R}}
\newcommand{\T}{\hbox{\bf t}}
\newcommand{\N}{\hbox{\bf n}}
\newcommand{\B}{\hbox{\bf b}}
\newcommand{\C}{\hbox{\bf c}}
\newcommand{\X}{\hbox{\bf X}}
\newtheorem{theorem}{Theorem}[section]
\newtheorem{corollary}[theorem]{Corollary}
\newtheorem{definition}[theorem]{Definition}
\newtheorem{remark}[theorem]{Remark}
\begin{document}

\title{Special Weingarten surfaces foliated by circles\footnote{Partially
supported by MEC-FEDER 
 grant no. MTM2004-00109.}}
\author{Rafael L\'opez\\
Departmento de Geometría y Topología\\
Universidad de Granada\\
18071 Granada (Spain)\\
e-mail:{\tt rcamino@ugr.es}}
\date{}
\maketitle


\begin{abstract} In this paper we study surfaces 
in Euclidean 3-space foliated by  pieces of circles and 
that satisfy a Weingarten condition of type 
$a H+b K=c$, where $a,b$ and $c$ are constant and  $H$ and $K$ denote 
the mean curvature and the Gauss curvature respectively. We prove that a such surface must be 
a surface of revolution, a  Riemann minimal surface or a generalized cone.
\end{abstract}

%%%%%%%%%%%%%%%%%%%%%%%%%%%%%%%%%%%%%%%%%%%%%
\section{Introduction}\label{intro}
%%%%%%%%%%%%%%%%%%%%%%%%%%%%%%%%%%%%%%%%%%%%
 
A surface $S$ in 3-dimensional Euclidean space $\r^3$ is called a {\it 
Weingarten surface} if there is some relation between its
two principal curvatures $\kappa_1$ and $\kappa_2$, 
that is, if there is a smooth function $W$ of two variable such that
$W(\kappa_1,\kappa_2)=0$. In particular, if $K$ and $H$  denote respectively  the Gauss and the mean curvature of  $S$, $W(\kappa_1,\kappa_2)=0$ implies a relation $U(K,H)=0$.  
In this work we study Weingarten surfaces that satisfy the simplest case for $U$, that is, that $U$ is of linear type:
\begin{equation}\label{w1}
a H+b K=c,
\end{equation}
where $a, b$ and $c$ are constant with $a^2+b^2\not=0$. We say then that $S$ is a {\it special Weingarten surface} and we abbreviate it by SW-surface. Constant mean curvature 
surfaces  ($b=0$) or constant Gauss curvature surfaces ($a=0$) are SW-surfaces.  The classification of Weingarten surfaces is almost completely open today. Weingarten 
introduced this kind of surfaces in the context of the problem of finding all surfaces isometric to a given surface of revolution \cite{we1,we2}. Along the history,
 they  have been of interest for geometers 
\cite{ch,hw,ho,ho2,vo} and more recently in \cite{gmm,ks,lo3,pa,rs}. Applications of Weingarten surfaces  on computer aided design and  shape investigation  can seen in \cite{bg}.

Among all  surfaces, a first class is the one of  surfaces of revolution. In such case, 
Equation (\ref{w1}) leads to an ordinary differential equation and its study is then simplified to find the profile curve (for example, see \cite{pa} if $c=0$).  
Recall that in the case of constancy of 
mean curvature or Gauss curvature, they are well known \cite{de,ei}. 
A more general family of surfaces of revolution are the cyclic surfaces. 

\begin{definition} A cyclic surface in Euclidean space $\r^3$ is a surface determined by a smooth uniparametric family of pieces of circles.
\end{definition}

Thus, a cyclic surface is a surface  foliated by pieces of circles meaning
that there is a one-parameter family of planes which meet $S$  in these 
pieces of circles.  The planes are not assumed parallel, and
if two circles should lie in planes that happen to be parallel, the
circles are not assumed coaxial.  We  point out that a sphere is a surface 
 such that any family of planes
(parallel or not) intersects it  in circles.
 The study of cyclic hypersurfaces with constant curvature in different ambient spaces was re-opened recently by  Jagy \cite{ja1,ja2}. See also \cite{lo2}.

The aim of this paper is if, besides the surfaces of revolution, there exist 
 new cyclic SW-surfaces. Our work is motivated by it happens in  the cases of 
constant mean curvature or constant Gauss curvature and that can summarized as follows. In both settings, the surface is an open of a sphere or the planes of the foliation are parallel \cite{en,lo1,ni}. When the planes of the foliation are parallel, then either it is a 
subset of a surface of revolution or it is a subset of one of the following non-rotational 
surfaces:

\begin{enumerate} 
\item It is one of the examples of minimal surfaces discovered by Riemann \cite{ri}. This happens when $H=0$.

\item It is a generalized cone, that is, a cyclic surface where the  circle centres lie in a 
straight-line and the radius function is linear. In this case, $K=0$. Locally, it can be 
parametrized by
$\X(u,v)=(f(u),g(u),u)+r(u)(\cos{(v)},\sin{(v)},0)$, where $f, g$ and $r>0$ are 
linear functions on $u$ \cite{lo1}. 
\end{enumerate}

 The first result that we shall obtain here states that in a cyclic 
SW-surface, the foliation planes must be parallel (except the trivial case of a sphere).

\begin{theorem}\label{t1} If $S$ is a SW-surface foliated  by pieces of
circles lying in a one-parameter family of planes, then either $S$ is a
subset of a round sphere or the planes in the family are parallel.
\end{theorem}

Once proved this result, we continue the search of cyclic SW-surfaces 
 in the situation of parallel planes. The conclusion that we shall obtain is that the 
circles of the foliation must be coaxial, unless the known cases of constant 
mean curvature or constant Gauss curvature.

\begin{theorem} \label{t2} Let $S$ be a SW-surface foliated by 
pieces of circles lying in a one-parameter family of parallel planes. Then 
either $S$ is a piece of a surface of revolution or $S$ is part of one of the Riemann minimal examples or a generalized cone. 
\end{theorem}

As conclusion of Theorems \ref{t1} and \ref{t2}, we have

\begin{corollary} \label{co1} A surface in 
Euclidean space foliated by pieces of circles and that 
satisfies a condition of type $a H+bK=c$, where $a,b$ and $b$ are constant must be a surface of revolution, a Riemann minimal surface or a generalized cone.
\end{corollary}

Therefore, although the family of surfaces satisfying the equation 
$aH+bK=c$ is larger than the one of constant mean curvature and 
constant Gauss curvature and so,  one could think the  existence of cyclic non-rotational surfaces for each three real numbers 
$(a,b,c)$, Corollary \ref{co1} says that this only occurs for 
the known cases of $H=0$ or $K=0$. In this sense,  we can view these two class of surfaces as  a special set of surfaces in the family of SW-surfaces.
 
 \begin{corollary} Riemann examples of minimal surfaces  and generalized cones are the only non-rotational cyclic surfaces that satisfy a  Weingarten relation of type $aH+bK=c$. 
\end{corollary}

This paper is organized as follows. In Section 2, we recall some concepts of the 
classical differential geometry of surfaces in $\r^3$. Section 3 and 4 are devoted to show 
Theorems \ref{t1} and \ref{t2} respectively.

%%%%%%%%%%%%%%%%%%%%%%%%%%%%%%%%%%%%%%%%%%%%%%%
\section{Preliminaries}
%%%%%%%%%%%%%%%%%%%%%%%%%%%%%%%%%%%%%%%%%%%%%%%%

In this section, we fix some notation on local classical differential geometry on 
surfaces. Let $S$ be a surface in $\r^3$ and consider  $\X=\X(u,v)$ 
 a local parametrization of $S$ defined in the $(u,v)$-domain. 
Let $N$  denote the unit normal vector field on $S$ given by
$$N=\frac{\X_u\wedge \X_{v}}{|\X_{u}\wedge \X_{v}|},
\qquad \X_u=\frac{\partial \X}
{\partial u},\   \X_v=\frac{\partial \X}{\partial v},$$
where $\wedge$ stands the cross product of $\r^3$. 
In each tangent plane, the induced metric 
$\langle,\rangle$  is determined by the first
fundamental form
$$I=\langle d\X,d\X\rangle=Edu^2+2Fdudv+Gdv^2,$$
with differentiable coefficients
$$E=\langle \X_u,\X_u\rangle,\quad F=\langle\X_u,\X_v\rangle,\quad
G=\langle \X_v,\X_v\rangle.$$
The shape operator of the immersion is represented by the second fundamental
form
$$II=-\langle d N ,d\X\rangle=e\ du^2+2f\ dudv+g\ dv^2,$$
with 
$$e=\langle  N ,\X_{uu}\rangle,\quad f=\langle N ,\X_{uv}\rangle,
\quad g=\langle  N ,\X_{vv}\rangle.$$
Under this parametrization $\X$, the mean curvature $H$ and the Gauss curvature $K$ have the classical expressions
$$H=\frac{eG-2f f+gE}{2 (EG-F^2)},\hspace*{1cm}K=\frac{eg-f^2}{EG-F^2}.$$
Let us denote by $[, ,]$ the determinant in $\r^3$ and put
$W=EG-F^2$. Then $H$ and $K$ write  as  
\begin{equation}\label{mean}
H=\frac{G[\X_u,\X_v,\X_{uu}]-2F [\X_u,\X_v,\X_{uv}]+E [\X_u,\X_v,\X_{vv}]}{2W^{3/2}}:=\frac{H_1}{2W^{3/2}},
\end{equation}
\begin{equation}\label{gauss}
 K=\frac{[\X_u,\X_v,\X_{uu}] [\X_u,\X_v,\X_{vv}] -[\X_u,\X_v,\X_{uv}]^2}{W^2}:=\frac{K_1}{W^2}.
\end{equation}
Using (\ref{mean}) and (\ref{gauss}), a SW-surface satisfies the condition
$$a\frac{H_1}{2W^{3/2}}+b\frac{K_1}{W^2}=c$$
or, equivalently, 
$$aH_1 W^{1/2}=2(c W^2-b K_1).$$
Squaring both sides, we have
\begin{equation}\label{w3}
a^2 H_1^2 W-4(c W^2-b K_1)^2=0.
\end{equation}

The proof of Theorems \ref{t1} and \ref{t2} involves  explicit computations 
of identity (\ref{w3}) and subsequent manipulations. As we shall see in the 
next two sections, Equation 
(\ref{w3}) reduces to an expression that   can be written as a linear combination of the functions 
$\cos{(jv})$ and $\sin{(j v)}$, $0\leq j\leq 8$, whose coefficients $A_j$ and $B_j$ 
 are functions on the $u$-variable. Therefore, they must vanish in some $u$-interval.
The work then is to compute explicitly these coefficients by successive 
manipulations. The author 
 was able to obtain the results using the  symbolic  program 
Mathematica to check his work. The computer was used in each calculation several times, giving understandable expressions of the coefficients $A_j$ and $B_j$. 

%%%%%%%%%%%%%%%%%%%%%%%%%%%%%%%%%%%%%%%%%%%%%%%
\section{Proof of Theorem \ref{t1}}
%%%%%%%%%%%%%%%%%%%%%%%%%%%%%%%%%%%%%%%%%%%%%%%%

In this section, we  follow the same ideas as in \cite{ni} and \cite{ja1} for the case that the mean curvature is constant. For this, we 
wish to construct an appropriate  coordinate system to our foliation of the 
surface $S$. Let us 
denote by $\Pi_u$  these planes in such way $S\cap \Pi_u$ is 
each piece of the circles of the foliation. Consider 
a smooth unit  vector field $Z$ that is normal to  the planes $\Pi_u$. 
Next, we  take a particular integral curve $\Gamma=\Gamma(u)$ of $Z$ parametrized by arclength, that is, 
$\T(u):=\Gamma'(u)=Z(\Gamma(u))$, where $\T$ is the unit tangent vector to 
$\Gamma$.   Consider the Frenet frame
of the curve $\Gamma$, $\{\T,\N,\B\}$, where $\N$ and $\B$ denote the normal and binormal vectors respectively. 

Locally we parametrize $S$ by
\begin{equation}\label{para}
\X(u,v)=\C(u)+r(u)(\cos{v}\ \N (u)+\sin{v}\ \B (u)),
\end{equation}
where $r=r(u)>0$ and $\C=\C(u)$ denote respectively 
the radius and centre of each $u$-circle of the foliation. Consider the Frenet equations of the curve $\Gamma$:
\begin{eqnarray*}
\T'&=&\hspace*{.8cm}\kappa \N\\
\N'&=&-\kappa \T+\ \ \sigma \B\\
\B'&=&\hspace*{.6cm}-\sigma\N
\end{eqnarray*}
where the prime $'$ denotes the derivative with respect to the $u$-parameter and $\kappa$ and $\sigma$ are the curvature and torsion of $\Gamma$, respectively.

Also, set
\begin{equation}\label{alfa}
\C'=\alpha\T+\beta\N+\gamma\B,
\end{equation}
where $\alpha,\beta,\gamma$ are smooth functions on $u$.

By using the Frenet equations and (\ref{alfa}), a straightforward computation of 
(\ref{w3}) shows that it  can be expressed by a trigonometric polynomial on 
$\cos{(jv)},\sin{(jv)}$. Exactly, there 
exist smooth functions on $u$, namely $A_j$ and $B_j$, such that Equation (\ref{w3}) writes as
\begin{equation}\label{formula}
A_0+\sum_{j=1}^8 \Bigl(A_j(u) \cos{(jv)} +B_j(u)\sin{(j v)}\Bigr)=
0,\hspace*{.5cm}u\in I, v\in J
\end{equation}
Since this is an expression on the independent trigonometric terms $\cos{ (j v)}$ and $\sin{(jv)}$, all coefficients $A_j,B_j$ must be zero.

In the reasoning to prove Theorem \ref{t1}, we shall assume 
that the planes of the foliation are not parallel and then, our objective will be to show 
that the surface is included in a sphere of Euclidean space. Thus, in our assumption, the 
integral curve $\Gamma$ is not a straight line perpendicular to each $u$-plane and so, 
the curvature $\kappa$  is not vanishing. The surface $S$ is part of a sphere if and only if each point of $S$ is the same distance from a fixed point $c_0$ in $\r^3$. We shall recognize  it in a few number of situations as the following: (i) the centre curve $\C$ is constant as well as the radius function $r(u)$ or; (ii) the curve $\C$ can written as 
a combination of $\T,\N$ and $\B$ in such way that the parametrization (\ref{para}) is now 
$\X(u,v)=c_0+\varphi(u,v)\T(u)+\phi(u,v)\N (u)+\psi(u,v)\B(u)$ where 
$|\X(u,v)-c_0|=\sqrt{\varphi^2+\phi^2+\psi^2}$ is a non-zero constant, the radius of the sphere that we are looking for.

In the proof, we distinguish two situations according the value $c$ in  (\ref{w1}).

%%%%%%%%%%%%%%%%%%%%%%%%%%%%%%%%%%%%%%%%%%%%%%%
\subsection{Case $c=0$ in the relation $aH+bK=c$.}
%%%%%%%%%%%%%%%%%%%%%%%%%%%%%%%%%%%%%%%%%%%%%%%%

Without loss of generality we assume that $a=1$. The coefficients $A_8$ and $B_8$ are 
\begin{eqnarray}\label{a8}
A_8&=&\frac{1}{32}\kappa^2 r^8(\beta^6-(15\gamma^2+\kappa^2(b^2-3r^2))\beta^4\nonumber\\
&-&(15\gamma^4+6\gamma^2\kappa^2(b^2-3r^2)+\kappa^4 r^2(-2b^2+3r^2))\beta^2\nonumber\\
&+&(\gamma^2-\kappa^2 r^2)(\gamma^2+\kappa^2(b^2-r^2)).
\end{eqnarray}
\begin{eqnarray}\label{b8}
B_8&=&\frac{1}{16}\beta\gamma\kappa^2 r^8(3 \beta^4-2\beta^2(5\gamma^2+\kappa^2 (b^2-3 r^2))\nonumber\\
&+ &(\gamma^2-\kappa^2 r^2)(3\gamma^2+\kappa^2(2 b^2-3r^2))),
\end{eqnarray}
From $B_8=0$, we discuss three cases.
\begin{enumerate}
\item Case $\beta=0$ in some sub-interval of $I$. Then
$$A_8=-\frac{1}{32}\kappa^2 r^8(\gamma^2-\kappa^2 r^2)^2(\gamma^2+\kappa^2(b^2-r^2)).$$
If $\gamma^2= \kappa^2 r^2$, then
$$A_6=-\frac98 b^2\kappa^6 r^{10}(\alpha^2-r'^2),\hspace*{.5cm}B_6=\pm\frac94 b^2 \alpha\kappa^6 r^{10} r',$$
which it is implies that $\alpha=0$ and $r$ is a constant function. But then $A_4=-2 r^{12} b^2\kappa^8$, obtaining a contradiction. As conclusion, and from $A_8=0$, we 
have that $\gamma=\pm\kappa\sqrt{r^2-b^2}$. Now 
$$A_7=-\frac{1}{16}b^4\alpha\kappa^7 r^9,\hspace*{.5cm}B_7=\pm\frac{b^4 \kappa^7 r^{10} r'}{16\sqrt{r^2-b^2}}.$$
Thus $r$ is a constant function and $\alpha=0$. Then 
$A_5=\pm r^9\kappa^7 b^4\tau\sqrt{r^2-b^2}/4$, that is, $\tau=0$ or $r=\pm b$. If $\tau=0$,  $A_4=r^8b^2(r^2-b^2)(5r^2+3b^2)\kappa^8/8$, which it leads to $r=\pm b$, and so, $\gamma=0$. Then 
(\ref{alfa}) implies that the curve of centres $\C(u)$ is constant, $\C(u)=c_0$ for 
some $c_0\in\r^3$. From (\ref{para}), 
$|\X(u,v)-c_0|=b^2$, that is, our surface is an open of a sphere of radius $|b|$. 

We  summarize this case by saying that if the foliation planes are not parallel, then the surface is a piece of a sphere.

\item Case $\gamma=0$ in some sub-interval of $I$. The coefficient $A_8$ is 
$$A_8=\frac{1}{32}\kappa^2 r^8(\beta^2+\kappa^2 r^2)^2(\beta^2+\kappa^2(r^2-b^2)).$$
Then $\beta^2=\kappa^2(b^2-r^2)$. Without loss of generality, we assume that 
$\beta=\kappa\sqrt{b^2-r^2}$. It follows that
$$A_7=-\frac{1}{16} b^4 \kappa^7 r^9\left(\alpha+\frac{r r'}{\sqrt{b^2-r^2}}\right).$$
Then $\alpha= -rr'/\sqrt{b^2-r^2}$. From (\ref{alfa}), we can  write
$$\C'=-\frac{r r'}{\sqrt{b^2-r^2}}\T+\kappa\sqrt{b^2-r^2}\N=
(\sqrt{b^2-r^2}\T)'.$$
Then there exists $c_0\in\r^3$ such that 
$\C=c_0+\sqrt{b^2-r^2}\T$ and the parametrization $\X$ of the surface is now
$$\X(u,v)=c_0+\sqrt{b^2-r^2}\T+r(\cos({v)}\N+\sin{(v)}\B).$$
This implies that $|\X(u,v)-c_0|=b^2$ and $S$ is again a piece of 
a sphere of radius $|b|$. In this setting, the same conclusion is obtained 
as in the above case.

\item Case $\beta\gamma\not=0$. From $B_8=0$ in (\ref{b8}), we can calculate $\beta^2$:
\begin{equation}\label{beta}
\beta^2=\frac13(5\gamma^2+b^2\kappa^2-3\kappa^2 r^2\pm A),
\end{equation}
where $A=\sqrt{16\gamma^4+4b^2\gamma^2\kappa^2+ b^4\kappa^4-12\gamma^2\kappa^2 r^2}$.
We consider the sign '+' in the value of $\beta^2$ (similarly with the choice $-$). 
Let us put it into $A_8$ and taking into account that $\kappa\not=0$, we obtain the following identity:
$$416\gamma^6+b^6\kappa^6+96\gamma^4\kappa^2(b^2-3r^2)+18\gamma^2\kappa^4 b^4=-(b^4\gamma^4+112\gamma^4+16\gamma^2\kappa^2(b^2-3r^2))A.$$
Squaring both sides and after some manipulations, we obtain
$$(\gamma^2-\kappa^2 r^2)\Bigl((4\gamma^2+b^2\kappa^2)^2-16\gamma^2\kappa^2 r^2\Bigr)=0.$$
We discuss each one of the possibilities:

\begin{enumerate}
\item $\gamma^2=\kappa^2 r^2$.  Using (\ref{beta}),  $\beta^2=2\kappa^2(b^2+2r^2)/3$ and returning with the 
computation of $A_8$ in (\ref{a8}), we have
$$A_8=-\frac{1}{216}\kappa^8 r^8(b^2+2r^2)(b^2+8r^2)^2.$$
Then $A_8=0$ yields a contradiction.
\item $(4\gamma^2+b^2\kappa^2)^2-16\gamma^2\kappa^2 r^2=0$. From here, we obtain the value of 
$\gamma^2$:
$$\gamma^2=\frac{\kappa^2}{4}\left(r\pm\sqrt{r^2-b^2}\right)^2.$$
Then the value of $\beta^2$ in (\ref{beta}) is
$$\beta^2=\frac{1}{12}\kappa^2(-5 b^2+2r^2+14 r\sqrt{r^2-b^2}).$$
From (\ref{a8}), Equation $A_8=0$ gives
$$(r^2-b^2)(b^4-14 b^2 r^2+16 r^4+(16 r^3-6r b^2)\sqrt{r^2-b^2})=0,$$
in particular, $r$ is a constant function.
The manipulation with the second factor implies that it cannot vanish. Thus 
$r^2=b^2$. But then $\beta^2=-\kappa^2 b^2/4$: contradiction. 
\end{enumerate}
\end{enumerate}

%%%%%%%%%%%%%%%%%%%%%%%%%%%%%%%%%%%%%%%%%%%%%%%
\subsection{Case $c\not=0$ in the relation $aH+bK=c$.}
%%%%%%%%%%%%%%%%%%%%%%%%%%%%%%%%%%%%%%%%%%%%%%%%

 Without loss of generality, we shall assume that $c=1$. The computation
of the coefficients $A_8$ and $B_8$ gives
\begin{equation}
A_8=-\frac{1}{32}r^8 x_1,\hspace*{1cm}B_8=\frac{1}{16}\beta\gamma r^8 x_2,
\end{equation}
where
\begin{eqnarray}\label{x1}
x_1&=&\beta^8-(28\gamma^2+\kappa^2(a^2+2b-4r^2))\beta^6\nonumber\\
&+&(70\gamma^4+15\gamma^2\kappa^2(a^2+2b-4r^2)+\kappa^4(b^3-3(a^2+2b)r^2+6r^4))\beta^4\nonumber\\
&+&(-28\gamma^6-15\gamma^4\kappa^2(a^2+2b-4r^2)+\kappa^6 r^2(2b^2-3(a^2+2b)r^2+4r^4)\nonumber\\
&-&
6\gamma^2\kappa^4(b^2-3(a^2+2b)r^2+6r^4))\beta^2\nonumber\\
&+&(\gamma^2-\kappa^2 r^2)^2(\gamma^4+\gamma^2\kappa^2(a^2+2b-2r^2)+\kappa^4(b^2-(a^2+2b)r^2+r^4)).
\end{eqnarray}
\begin{eqnarray}\label{x2}
x_2&=&-4\beta^6+(28\gamma^2+3\kappa^2(a^2+2b-4r^2)\beta^4\nonumber\\
&-&2 
(14\gamma^4+5\gamma^2\kappa^2(a^2+2b-4r^2)+\kappa^4(b^2-3(a^2+2b)r^2+6r^4)\beta^2\nonumber\\
&+&(\gamma^2-\kappa^2 r^2)(4\gamma^4+\gamma^2\kappa^2(3a^2+6b-8r^2)+\kappa^4(2b^2-
3(a^2+2b)r^2+4r^4).
\end{eqnarray}
We discard the cases  $a=0$ or $b=0$, corresponding to 
the known situations of (non-zero) constant mean or Gauss curvature: in such 
case,   $S$ is a piece of a sphere. 
From $B_8=0$, we discuss the following cases:
\begin{enumerate}
\item Case $\gamma=0$. From $A_8=0$, 
$$\beta^4-\kappa^2(a^2+2b-2r^2)\beta^2+\kappa^4(b^2-(a^2+2b)r^2+r^4))=0.$$
Then 
$$\beta^2=\frac12\kappa^2\left(a^2+2b-2r^2\pm a\sqrt{a^2+4b}\right).$$
In particular, $a^2+4b\geq 0$. In the reasoning, we shall suppose  the positive sign in $\pm$ of the expression of $\beta^2$. According to the value of $a^2+4b$, we distinguish two cases.

\begin{enumerate}
\item $a^2+4b=0$. From the value of $\beta^2$, the coefficient $B_5$ gives
$$B_5=\frac{1}{128}a^4\kappa^5 r^7\sqrt{a^2-4r^2}\left(\alpha\sqrt{a^2-4r^2}+2 r r'\right)^2=0.$$
If $\sqrt{a^2-4r^2}=0$, $\beta=0$ and we  are in the case "$\beta=0$" that it 
will be studied in the second case of this subsection. Thus,  
$$\alpha=-\frac{2 r r'}{\sqrt{a^2-4 r^2}}$$
and this allows us to write
$$\C'=(\frac{\sqrt{a^2-4r^2}}{2}\T)'.$$
As a consequence, there exists some fixed vector $c_0$ such that 
$$\C=c_0+\frac{\sqrt{a^2-4r^2}}{2}\T.$$
Then the parametrization $\X$ in (\ref{para}) gives $|\X(u,v)-c_0|^2=a^2/4$, that is, the surface is a open of a sphere of radius $|a|/2$.

In the second case, the computation of the coefficient $B_5=0$ implies
$$\alpha\sqrt{a^2-4r^2}+2 r r'=0\hspace*{.5cm}\mbox{or}\hspace*{.5cm}\tau=0.$$
obtaining the same result as above.
\item $a^2+4b>0$. The coefficient $A_7$ is
$$B_7=\frac{1}{64}aAB\kappa^5 r^9 (\alpha\kappa^2-\kappa\beta'+\kappa'\beta)=0,$$
with 
$$A=2b+a(a+\sqrt{a^2+4b})\hspace*{1cm}B=a^3+4ab+(a^2+2b)\sqrt{a^2+4b}).$$
The number $A$ does not vanish and $B=0$ only if $a^2+4b=0$. As 
conclusion,
$$\alpha=\left(\frac{\beta}{\kappa}\right)'.$$
Following (\ref{alfa}),  the derivative of the curve $\C$ is 
$$\C'=\left(\frac{\beta}{\kappa}\right)'\T+\beta\N=\left(\frac{\beta}{\kappa}\T\right)'.$$
From (\ref{para}), the parametrization of the surface is 
$$\X(u,v)=c_0+\frac{\beta}{\kappa}\T+r(\cos{(v)}\N+\sin{(v)}\B),$$
for some fixed vector $c_0$. Using the value of $\beta^2$, we have, 
$$|\X(u,v)-c_0|^2=\frac{\beta^2}{\kappa^2}+r^2=\frac12 a^2+b+\frac{a}{2}\sqrt{a^2+4b}.$$
This means that the surface is an open of a certain sphere.
\end{enumerate}

\item Case $\beta=0$. Now 
$$A_8=-\frac{1}{32}r^8(\gamma^2-\kappa^2 r^2)^2 x_1,\hspace*{1cm}
A_7=-\frac{1}{16}\alpha\kappa r^9(\gamma^2-\kappa^2 r^2) y_1,$$
where
$$x_1=\gamma^4+\kappa^2\gamma^2(a^2+2b-2r^2)+
\kappa^4(b^2-(a^2+2b)r^2+r^4).$$
$$y_1=8\gamma^4+(7(a^2+2b)-16 r^2)\kappa^2\gamma^2+\kappa^4(6 b^2-7(a^2+2b)r^2+8r^4).$$
We discuss three possibilities:
\begin{enumerate}
\item Case $\gamma^2=\kappa^2 r^2$. Then
$$A_6=\pm\frac98 b^2\kappa^6 r^{10}(\alpha^2-r'^2),\hspace*{1cm}B_6=\pm\frac98 b^2\alpha\kappa^6 r^{10} r'.$$
Then $\alpha=0$ and $r$ is a constant function. Then $A_4=-2r^{12}b^2\kappa^8$, giving 
a contradiction.
\item Case $x_1=\alpha=0$. We know that 
$$\gamma^2=\frac12 k^2(\pm a\sqrt{a^2+4b}-(a^2+2b-2r^2)).$$
Without loss of generality, we assume the sign + in $\pm$. Then 
$$B_7=\frac{1}{32}a r^9 C\kappa^3(\gamma^2-\kappa^2 r^2)(\kappa\gamma'-\kappa'\gamma),$$
where $C=a^3+4ab-(a^2+2b)\sqrt{a^2+4b}$. For each pair $(a,b)$ of real numbers, 
$C\not=0$ except when $a^2+4b=0$. Then $\kappa\gamma'-\kappa'\gamma=0$. From this 
equation, we conclude that $r$ is a constant function. We discuss both situations according to the value of $a^2+4b$:
\begin{enumerate}
\item Let $a^2+4b=0$. Equation $A_4=0$ implies $16 r^4+8 a^2r^2-3a^4=0$, that is, $4r^2-a^2=0$. Thus $\gamma=0$ and this 
case was studied above.

\item Let $a^2+4b>0$. Now $A_5=0$ implies $\tau=0$. After some manipulations, $A_4=0$ and $A_2=0$ give $\kappa=0$:  contradiction.
\end{enumerate}

\item Case $x_1=y_1=0$. Then $8x_1-y_1=0$ means
$$(a^2+2b)\gamma^2+\kappa^2(2b^2-(a^2+2b) r^2)=0.$$
In particular, $a^2+2b\not=0$ and $\gamma^2=\kappa^2((a^2+2b)r^2-2b^2)/(a^2+2b)$. With this 
value of $\gamma^2$, 
$$x_1=-\frac{a^2 b^2(a^2+4b)\kappa^4}{(a^2+2b)^2}=0,$$
which it is a contradiction.
\end{enumerate}

\item Case $\beta\gamma\not=0$. This case is more difficult in the computations due to that
the expressions are very cumbersome. We only give the details. It follows from the
expressions of $A_8$ and $B_8$ in (\ref{a8}) that 
$x_1=x_2=0$: see (\ref{x1}) and (\ref{x2}). We begin to compute the value 
of $\beta^2$. For this, we define $x_3:=4x_1+\beta^2 x_2=0$ 
and $x_4:=4x_3-(84\gamma^2+\kappa^2(a^2+2b-4r^2))x_2$.  Now $x_4$ is a $2$-degree 
polynomial on $\beta^2$ and we can calculate $\beta^2$:
\begin{equation}\label{beta2}
\beta^2=\frac{\xi\pm\sqrt{\xi^2-\eta \lambda \zeta}}{\eta},
\end{equation}
where
\begin{eqnarray*}
\xi&=&960\gamma^6+320\gamma^4\kappa^2(a^2+2b-4r^2)+\kappa^6(b^2(a^2+2b)-(3a^4+12 a^2 b+4 b^2)r^2)\\
&+&5\gamma^2\kappa^4(a^4+4 a^2 b+12 b^2-32 r^2(a^2+2b-2r^2)),\\
\eta&=&1344\gamma^4+(3a^4+12a^2 b+4a^2)\kappa^4+80\gamma^2\kappa^2(a^2+2b-4r^2),\\
\lambda&=&\gamma^2-\kappa^2 r^2,\\
\zeta&=&320\gamma^6+80\gamma^4\kappa^2(3a^2+6b-6r^2)+\kappa^6(2b^2(a^2+2b)-(3a^4+12a^2 b+4b^2) r^2)\\
&+&\gamma^2\kappa^4(3a^4+12a^2 b+164 b^2-240(a^2+2b) r^2+320 r^4).
\end{eqnarray*}
For each one of the two values of $\beta^2$, we return to $x_2=0$ obtaining the following:
\begin{equation}\label{betagamma}
\lambda\mu \eta^3  \rho =0,
\end{equation}
where
\begin{eqnarray*}\mu&=&16\gamma^4+a^2(a^2+4b)\kappa^4+8\gamma^2\kappa^2(a^2+2b-2r^2),\\
\rho&=&(16\gamma^4+32\kappa r\gamma^3+4\kappa^2(a^2+2b+4r^2)\gamma^2+4(a^2+2b)\kappa^3 r\gamma+b^2\kappa^4)\\
& &(16\gamma^4-32\kappa r\gamma^3+4\kappa^2(a^2+2b+4r^2)\gamma^2-4(a^2+2b)\kappa^3 r\gamma+b^2\kappa^4).
\end{eqnarray*}
From Equation (\ref{betagamma}), we have four cases to discuss in such way that we can compute the value of $\gamma^2$ and, next, putting it  in (\ref{beta2}), the value  of $\beta^2$. For instance, if $\lambda=0$, that is,  $\gamma^2-\kappa^2 r^2=0$,
 the value of
 $\beta^2$    is
$$\beta^2=2\kappa^2\frac{2r^2(a^4+2a^2 b+28 b^2+80(a^2+2b) r^2)+b^2(a^2+2b)}{3a^4+12 a^2+b+4b^2+80(a^2+2b) r^2+1024 r^4}.$$
On the other hand,  $x_2$ writes now as
$$x_2=\beta^2(-4\beta^4+\beta^2\kappa^2(3a^2+6b+16 r^2)-2\kappa^4(b^2+2(a^2+2b)r^2))$$
and it follows that 
$$\beta^2=\frac{\kappa^2}{8}(3a^2+6b+16r^2\pm\sqrt{9a^4+36 a^2 b+4b^2+32 a^2 r^2+64 b r^2+256 r^4}).$$
Comparing both values of $\beta^2$, we know then $r^2$ is one of the following values:
$$r^2=-\frac{a^2(a^2+4b)}{8(a^2+2b)},\hspace*{.5cm}r^2=\frac{-1}{16}(a^2+2b\pm\sqrt{a^2(a^2+4b)}).$$
If we analyze, for example,  the first value of 
$r^2$, we know that 
$$\beta^2=\frac12(a^2+2b)\kappa^2\hspace*{.5cm}
\gamma^2=-\frac{a^2(a^2+4b)}{8(a^2+2b)}\kappa^2.$$
Now $A_7=0$ implies
$\alpha\kappa^7=0$, that is,  
$\alpha=0$. Equations $A_5=0$ and $B_5=0$ give  $\tau\kappa^7=0$, and so, 
$\tau=0$. With $j=4$, 
$$A_4=(21 a^6+130 a^4 b+240 a^2 b^2+96 b^3)\kappa^8=0,$$
$$B_4=(21 a^4+88 a^2 b+96 b^2)\kappa^8=0,$$
which it would imply $\kappa=0$, obtaining the desired contradiction.

\end{enumerate}

\begin{remark} Throughout our reasoning in the case $c\not=0$, 
it has appeared, as a particular case, that  $a^2+4b>0$ (or $a^2+4bc>0$ in the Weingarten relation (\ref{w1})). This is not casual: Weingarten surfaces that satisfy this property were treated by the very Hopf in \cite{ho2} by their special properties (see also \cite{gmm}).
\end{remark}

%%%%%%%%%%%%%%%%%%%%%%%%%%%%%%%%%%%%%%%%%%%%%%%%%%%
\section{Proof of Theorem \ref{t2}}
%%%%%%%%%%%%%%%%%%%%%%%%%%%%%%%%%%%%%%%%%%%%%%%%%%%%

Once we have proved Theorem \ref{t1}, we consider SW-surfaces foliated 
by pieces of circles in parallel planes. The conclusion that we shall arrive is that 
either (i) the circles of the foliation must coaxial, that is, the surface is an open subset of a surface of revolution or (ii) it is part of a Riemann minimal surface ($H=0$) or of a generalized cone ($K=0$). This is the statement of Theorem \ref{t2}, which it will be proved in this section. Because our reasoning is of local character, we can assume the planes are parallel to the $x_1 x_2$-plane. Therefore, the surface $S$  writes as 
$$\X(u,v)=(f(u),g(u),u)+r(u)(\cos{v},\sin{v},0),$$
where  $f, g$ and $r$ are smooth functions in some $u$-interval $I$ and $r>0$ denotes the 
radius of each circle of the foliation. With this parametrization,  $S$ is a surface of 
revolution if and only if $f$ y $g$ are constant functions. 
If we compute (\ref{w3}), we obtain  
\begin{equation}\label{p0}
\sum_{j=0}^{8} A_j(u) \cos{(j v)}+B_j(u)\sin{(j v)}=0.
\end{equation}
Again, the functions $A_j$ and $B_j$ on $u$ must vanish on $I$. 

In our reasoning, we shall assume that the foliated circles are not coaxial and that 
that $b^2+c^2\not=0$ and $a^2+c^2\not=0$ (which it would yield that $S$ is part of a Riemann example or of a generalized cone). With these assumptions, we will arrive to a contradiction.
As in the above section, we distinguish two cases according to the value of $c$ in the relation 
$aH+bK=c$.

 %%%%%%%%%%%%%%%%%%%%%%%%%%%
\subsection{Case $c=0$.}
%%%%%%%%%%%%%%%%%%%%%%%%%%%

 In this particular situation,  the sum in (\ref{p0}) is until $j=4$. 

\begin{enumerate}
\item First, we  consider the cases that one of the functions $f$ or $g$ is  constant.
 For simplicity, we shall consider $f'=0$ in some interval. Then $A_4 $ writes as 
$$A_4=\frac18 a^2 r^6 g'^2(r g''-2 r' g').$$
As $g'\not=0$, we have that $r g''-2 r' g'=0$. Then 
$g'=\lambda r^2$ for some positive constant $\lambda\not=0$. Now
$$A_2=-\frac12\lambda^2 r^8(a^2 r^2 A^2-16 b^2 r'^2),\hspace*{.5cm}
B_1=2\lambda r^7 r'(a^2 r A^2-8b^2 r''),$$
where
$$A=1+\lambda^2 r^4+r'^2-r r''.$$
From Equation $A_2=0$ and the value of $A$, we discard the case that $r$ is a constant 
function. Thus, 
the combination of $A_2=0$ and $B_1=0$ leads to
that the function $r$ satisfies the ordinary differential equation $r r''-2 r'^2=0$. Then 
$$r(u)=\frac{\alpha}{u+\beta},\hspace*{1cm}\alpha,\beta\in\r.$$
But then $A_2=0$ gives a polynomial on $u$ given by
$$16b^2(u+\beta)^6-a^2((u+\beta)^4-\alpha^2+\lambda^2\alpha^4)^2=0,\ \ \ \forall u\in I$$
and whose leading coefficient, corresponding to $u^8$, is $-a^2$: this is a
 contradiction. This means that the assumption that $f$ is constant is impossible.

\item We  assume that both $f$ and $g$ are not constant functions. Then $f',g'\not=0$. The computation of 
$B_4$ gives now
\begin{eqnarray*}
B_4&=&\frac14 a^2 r^6\Bigl(rg'f''+f'(-4g'r'+r g'')\Bigr)\Bigl(-2 f'^2 r'+r f' f''+g'(2g' r'-r g'')\Bigr).
\end{eqnarray*}
We have two possibilities.
\begin{enumerate}
\item Case $rg'f''+f'(-4g'r'+r g'')=0$. Then
\begin{equation}\label{efe}
f''=\frac{f'(4 g' r'-r g'')}{r g'},
\end{equation}
and the coefficient $A_4=0$ gives
$$A_4=\frac{a^2 r^6 (f'^2+g'^2)^2}{8 g'^2}(r g''-2 g' r')^2.$$
Thus $r g''-2 g' r'=0$, that is, $g'=\lambda r^2$ with $\lambda>0$. Using (\ref{efe}),
 the same occurs for 
$f'$: $f'=\mu r^2$, $\mu>0$. The computation of $A_2$ and $A_1$ leads to $$A_2=-\frac12(\lambda^2-\mu^2)r^8(a^2 r^2 A^2-16 b^2 r'^2),$$
$$A_1=2\mu r^7 r' (-8b^2 r''+a^2 r A^2),$$
where the value of $A$ is now
$$A=1+(\lambda^2+\mu^2) r^4+r'^2-r r''.$$
From the expression of $A$ together $A_2=0$, we conclude that $r$ is not 
a constant function. By combining 
$A_2=0$ and $A_1=0$, we obtain $r r''-2 r'^2=0$ again. The contradiction is obtained as in the  case that $f$ is a constant function.
\item Case $-2 f'^2 r'+r f' f''+g'(2g' r'-r g'')=0$. From here, we obtain $f''$ and putting it into $A_4$, it gives
$$A_4=-\frac{a^2 r^6 (f'^2+g'^2)^2}{8 f'^2}(r g''-2 g' r')^2.$$
Then $r g''-2 g' r'=0$: we are now in the position of the above case. 
\end{enumerate}

\end{enumerate}

%%%%%%%%%%%%%%%%%%%%%%%%%%%
\subsection{Case $c\not=0$.}
%%%%%%%%%%%%%%%%%%%%%%%%%%%

The computation of   $A_{8}$ and $B_{8}$ give respectively:
$$A_8=-\frac{1}{32}c^2 r^8(f'^8-28 f'^6 g'^2+70 f'^4 g'^4-28 f'^2 g'^6+g'^8).$$
$$B_8=\frac14 c^2 r^8 f' g' (-f'^6+7 f'^4 g'^2-7 f'^2 g'^4+g'^6).$$
 Since  $\alpha(u)=(f(u),g(u))$ is not a constant planar curve, we reparametrize 
it by the arclength, that is, $(f(u),g(u))=(x(\phi(u),y(\phi(u))$, where
$$f'(u)=\phi'(u)\cos{(\phi(u))},\hspace*{.5cm}g'(u)=\phi'(u)\sin{(\phi(u))},
\hspace*{.5cm}\phi'^2=f'^2+g'^2.$$
With this change of variable, the functions $A_8$ and $B_8$ write now as:
$$A_8=-\frac{1}{32} c^2 r^8 \phi'^8\cos{(8\phi(u))}.$$
$$B_8=-\frac{1}{32} c^2 r^8 \phi'^8 \sin{(8\phi(u))}.$$
As $c\not=0$ and $r>0$,  we conclude that $\phi'=0$ on  some interval, that is, 
$\alpha$ is a constant curve, obtaining a contradiction.

This finishes the proof of Theorem \ref{t2}.

We end with a  comment when $c=0$. For the cases of constant mean curvature or constant  Gauss curvature, the same above computations give:
\begin{enumerate} 
\item Constant mean curvature ($b=0$). 
$$f''=\lambda r^2,\hspace*{.5cm}g''=\mu r^2\hspace*{.5cm}
1+(\lambda^2+\mu^2) r^4+r'^2-r r''=0,$$
which it gives the Riemann examples of minimal surfaces ($\lambda^2+\mu^2\not=0$) and 
the catenoid ($\lambda=\mu=0$).
\item Constant Gauss curvature ($a=0$). 
$$f''=g''=r''=0,$$
that is, the functions $f, g$ and $r$ are linear on $u$ and so, 
the surface is a generalized cone. 
\end{enumerate}

\end{document}